*Статья приурочена к 60-и летию Ю.Е. Нестерова*

# Численные методы поиска равновесного распределения потоков в модели Бэкмана и модели стабильной динамики

*А.В. Гасников (ИППИ РАН, ПреМоЛаб МФТИ),*
*П.Е. Двуреченский (ПреМоЛаб МФТИ, WIAS, ИППИ РАН)*
*Ю.В. Дорн (Центр исследований транспортной политики, Институт экономики транспорта и транспортной политики, НИУ ВШЭ, ПреМоЛаб МФТИ, ИППИ РАН)*
*Ю.В. Максимов (ИППИ РАН, ПреМоЛаб МФТИ)*

**Аннотация**

В работе рассматриваются две модели транспортного равновесия: модель Бэкмана (1955) и модель стабильной динамики (Нестеров–де Пальма, 1998). В статье описаны эффективные численные процедуры поиска равновесия в этих моделях. Для модели Бэкмана будет использован метод Франк–Вульфа, а для модели стабильной динамики используется переход к двойственной задаче. Эта задача решается методом зеркального спуска с евклидовой прокс-структурой с помощью "рандомизации суммы". В работе также приводится другой способ решения (сглаженной) двойственной задачи. Этот способ базируется на современных вариантах метода ускоренного блочно-покомпонентного спуска. Такие подходы, насколько нам известно, представляются новыми. Кроме того, даже при использовании классического метода Франк–Вульфа, мы исходим из современных результатов о его сходимости.

## 1. Введение

В недавних работах [1, 2], посвященных сведению поиска равновесного распределения транспортных потоков на сетях к решению задач выпуклой оптимизации, было поставлено несколько таких задач со специальной "сетевой" структурой. Это означает, что, скажем, расчет градиента (стохастического градиента) функционала сводится к поиску кратчайших путей в графе транспортной сети. Эта специфика задач с одной стороны говорит о том, что в реальных приложениях размерность задачи может быть колоссально большой. Это связано с тем, что число реально используемых путей даже в планарном графе, как правило, пропорционально кубу числа вершин, а число вершин в реальных приложениях обычно не меньше тысячи, отметим, что в худших случаях число путей может расти экспоненциально с ростом числа вершин. С другой стороны, такого рода задачи имеют хорошую геометрическую интерпретацию, что позволяет эффективно снижать их размерность.

В частности, в разделе 2 мы описываем метод Франк–Вульфа [3 – 5] поиска равновесного распределения потоков в модели Бэкмана, который на каждой итерации требует решения задачи минимизации линейной функции от потоков на ребрах в сети (число ребер порядка нескольких тысяч) на прямом произведении симплексов (симплексов столько, сколько корреспонденций). Для реальных транспортных сетей (тысяча вершин) получается задача минимизации линейной функции в пространстве размерности миллиард, поскольку она зависит от распределения потоков по путям, число которых порядка миллиарда. Ясно, что если смотреть на эту задачу формально с точки зрения оптимизации, то все сводится к полному перебору миллиарда вершин всех симплексов (причем проработка одной вершины – это расчет соответствующего скалярного произведения, то есть порядка нескольких тысяч умножений). К счастью, транспортная специфика задачи позволяет с помощью алгоритма Дейкстры и более современных подходов [6] (в том числе учитывающих “планарность” сети: A*, ALT, SHARC, Reach based routing, Highway hierarchies, Contraction hierarchies и т.п. – этому планируется посвятить отдельную работу) решать описанную задачу делая не более десятка миллионов операций (типа умножения двух чисел с плавающей запятой), что намного быстрее. Такого рода конструкции возникают не только в связи с сетевой спецификой задачи [7], но именно для ситуаций, когда в задаче имеется сетевая структура, возможность такой редукции наиболее естественна и типична.

В разделе 3 мы предлагаем другой способ поиска равновесия в модели Бэкмана и аналогичных моделях (в модели стабильной динамики, в промежуточных моделях). Для этого мы переходим (следуя Ю.Е. Нестерову) к двойственной задаче, в которой целевой функционал оказывается зависящим только от потоков по ребрам, а не от распределения потоков по путям. Таким образом, задача сводится к поиску равновесного распределения потоков по ребрам. При этом в ходе вычислений потоков на ребрах, мы попутно (без дополнительных затрат) вычисляем порождающие их потоки по путям. Также в виду транспортно-сетевой специфики появляется возможность содержательной интерпретации [1] (подобно интерпретации Л.В. Канторовичем цен в экономике [8]), возникающих двойственных множителей, которые в ряде приложений представляют независимый самостоятельный интерес (например, в задаче о тарифной политике грузоперевозок РЖД [2] двойственные множители – тарифы, которые и надо рассчитывать). Также сетевая структура задачи дает возможность не рассчитывать градиент целевой функции на каждой итерации заново, а пересчитывать его, используя градиент, полученный на предыдущей итерации.

Грубо говоря, найдя кратчайшие пути, и посчитав на их основе градиент, мы сделаем шаг по антиградиенту, немного изменив веса ребер. Ясно, что большая часть кратчайших путей при этом останется прежними, и можно специально организовать их пересчет, чтобы ускорить вычисления. Похожая философия используется в покомпонентных спусках и в современных подходах к задачам huge-scale оптимизации [9, 10]. Однако сетевая структура задачи требует переосмысления этой техники, рассчитанной изначально в основном на свойства разреженности матриц, возникающих в условии задачи.

Применимость современных вариантов рандомизированных покомпонентных спусков к поиску равновесия в модели стабильной динамики, записанной в новой специальной форме, предложенной недавно Ю.Е. Нестеровым, изучается в разделе 4.

В заключительном разделе 5 приводится краткий сравнительный анализ описанных в статье методов.

Настоящая статья представляет собой одну из первых попыток авторов сочетать современные эффективные численные методы выпуклой оптимизации с сетевой структурой задачи, на примере задач, пришедших из поиска равновесного распределения потоков в транспортных сетях и сетях грузовых перевозок РЖД [1, 2]. В последствие планируется опубликовать еще несколько статей на эту тему, представляющую на, наш взгляд, большой интерес.

## 2. Метод Франк–Вульфа поиска равновесия в модели Бэкмана

Опишем наиболее популярную на протяжении более чем полувека модель равновесного распределения потоков Бэкмана [1, 11 – 15]. Первая половина этого раздела во многом повторяет текст раздела 4 работы [1].

Пусть транспортная сеть города представлена ориентированным графом $\Gamma = (V, E)$, где $V$ – узлы сети (вершины), $E \subset V \times V$ – дуги сети (рёбра графа), $O \subseteq V$ – источники корреспонденций ($S = |O|$), $D \subseteq V$ – стоки. В современных моделях равновесного распределения потоков в крупном мегаполисе число узлов графа транспортной сети обычно выбирают порядка $n = |V| \sim 10^3 - 10^4$. Число ребер $|E|$ получается в три четыре раза больше. Пусть $W \subseteq \{w = (i, j) : i \in O, j \in D\}$ – множество корреспонденций, т.е. возможных пар «источник» – «сток»; $p = \{v_1, v_2, ..., v_m\}$ – путь из $v_1$ в $v_m$, если $(v_k, v_{k+1}) \in E$, $k = 1, ..., m-1$, $m > 1$; $P_w$ – множество путей, отвечающих корреспонденции $w \in W$, то есть если $w = (i, j)$, то $P_w$ – множество путей, начинающихся в вершине $i$ и заканчивающихся в $j$;

$P = \bigcup_{w\in W} P_w$ – совокупность всех путей в сети $\Gamma$ (число “разумных” маршрутов $|P|$, которые потенциально могут использоваться, обычно растет с ростом числа узлов сети не быстрее чем $O\left(n^3\right)$ [12 – 14]); $x_p$ [автомобилей/час] – величина потока по пути $p$, $x = \left\{x_p:\ p\in P\right\}$; $f_e$ [автомобилей/час] – величина потока по дуге $e$:

$$f_e\left(x\right) = \sum_{p\in P}\delta_{ep}x_p\text{, где } \delta_{ep} = \begin{cases}1, & e\in p\\ 0, & e\notin p\end{cases}\quad (f = \Theta x);$$

$\tau_e\left(f_e\right)$ – удельные затраты на проезд по дуге $e$. Как правило, предполагают, что это – (строго) возрастающие, гладкие функции от $f_e$. Точнее говоря, под $\tau_e\left(f_e\right)$ правильнее понимать представление пользователей транспортной сети об оценке собственных затрат (обычно временных в случае личного транспорта и комфортности пути (с учетом времени в пути) в случае общественного транспорта) при прохождении дуги $e$, если поток желающих проехать по этой дуге будет $f_e$.

Рассмотрим теперь $G_p\left(x\right)$ – затраты временные или финансовые на проезд по пути $p$. Естественно считать, что $G_p\left(x\right) = \sum_{e\in E}\tau_e\left(f_e\left(x\right)\right)\delta_{ep}$.

Пусть также известно, сколько перемещений в единицу времени $d_w$ осуществляется согласно корреспонденции $w\in W$. Тогда вектор $x$, характеризующий распределение потоков, должен лежать в допустимом множестве:

$$X = \left\{x\ge 0:\ \sum_{p\in P_w} x_p = d_w,\ w\in W\right\}.$$

Рассмотрим игру, в которой каждой корреспонденции $w\in W$ соответствует свой, достаточно большой ($d_w \gg 1$), набор однотипных “игроков”, осуществляющих передвижение согласно корреспонденции $w$. Чистыми стратегиями игрока служат пути, а выигрышем – величина $-G_p\left(x\right)$. Игрок “выбирает” путь следования $p\in P_w$, при этом, делая выбор, он пренебрегает тем, что от его выбора также “немного” зависят $|P_w|$ компонент вектора $x$ и, следовательно, сам выигрыш $-G_p\left(x\right)$. Можно показать (см., например, [1]), что отыскание равновесия Нэша–Вардропа $x^*\in X$ (макро описание равновесия) равносильно решению задачи нелинейной комплементарности (принцип Вардропа):

$$\textit{для любых } w\in W,\ p\in P_w \textit{ выполняется } x_p^*\cdot\left(G_p\left(x^*\right) - \min_{q\in P_w} G_q\left(x^*\right)\right) = 0.$$

Действительно допустим, что реализовалось какое-то другое равновесие $\tilde{x}^* \in X$ , которое не удовлетворяет этому условию. Покажем, что тогда найдется водитель, которому выгодно поменять свой маршрут следования. Действительно, тогда

$$\text{существуют такие } \tilde{w} \in W,\ \tilde{p} \in P_{\tilde{w}}, \text{ что } \tilde{x}^*_{\tilde{p}} \cdot \left( G_{\tilde{p}}\left(\tilde{x}^*\right) - \min_{q \in P_{\tilde{w}}} G_q\left(\tilde{x}^*\right)\right) > 0 .$$

Каждый водитель (множество таких водителей не пусто, так как $\tilde{x}^*_{\tilde{p}} > 0$), принадлежащий корреспонденции $\tilde{w} \in W$ , и использующий путь $\tilde{p} \in P_{\tilde{w}}$ , действует не разумно, поскольку существует такой путь $\tilde{q} \in P_{\tilde{w}},\ \tilde{q} \neq \tilde{p}$ , что $G_{\tilde{q}}\left(\tilde{x}^*\right) = \min_{q \in P_{\tilde{w}}} G_q\left(\tilde{x}^*\right)$ . Этот путь $\tilde{q}$ более выгоден, чем $\tilde{p}$ . Аналогично показывается, что при $x^* \in X$ никому из водителей уже не выгодно отклоняться от своих стратегий.

Условие равновесия может быть переписано следующим образом [11 – 15]

$$\text{для всех } x \in X \text{ выполняется } \left\langle G\left(x^*\right), x - x^* \right\rangle \geq 0 .$$

Рассматриваемая нами игра принадлежит к классу, так называемых, потенциальных игр [16, 17], поскольку $\partial G_p(x)/\partial x_q = \partial G_q(x)/\partial x_p$ . Существует такая функция

$$\Psi\left(f(x)\right) = \sum_{e \in E} \int_0^{f_e(x)} \tau_e(z)\,dz = \sum_{e \in E} \sigma_e\left(f_e(x)\right),$$

где $\sigma_e\left(f_e\right) = \int_0^{f_e(x)} \tau_e(z)\,dz$ , что $\partial \Psi(x)/\partial x_p = G_p(x)$ для любого $p \in P$ . Таким образом, $x^* \in X$ – равновесие Нэша–Вардропа в этой игре тогда и только тогда, когда оно доставляет минимум $\Psi\left(f(x)\right)$ на множестве $X$ .

**Теорема 1 [1, 12 – 15].** *Вектор $x^*$ будет равновесием Нэша–Вардропа тогда и только тогда, когда*

$$x \in \operatorname{Arg}\min_x \left[ \Psi\left(f(x)\right) = \sum_{e \in E} \sigma_e\left(f_e(x)\right) :\ f = \Theta x,\ x \in X \right].$$

*Если преобразование $G(\cdot)$ строго монотонное, то равновесие $x$ единственно. Если $\tau_e'(\cdot) > 0$, то равновесный вектор распределения потоков по ребрам $f$ – единственный* (*это еще не гарантирует единственность вектора распределения потоков по путям $x$ [15]*).

Итак, будем решать задачу ( $\Psi_*$ – оптимальное значение функционала)

$$\Psi(f)=\sum_{e\in E}\sigma_e(f_e)\to\min_{\substack{f=\Theta x\\ x\in X}}$$

методом условного градиента [18 – 21] (Франк–Вульфа).

**Начальная итерация**

*Положим* $\tilde{t}_e^0=\partial\Psi(0)/\partial f_e=\tau_e(0)$ *и рассмотрим задачу*

$$\sum_{e\in E}\tilde{t}_e^0 f_e\to\min_{\substack{f=\Theta x\\ x\in X}} .$$

*Эту задачу можно переписать, как*

$$\min_{x\in X}\sum_{e\in E}\tilde{t}_e^0\sum_{p\in P}\delta_{ep}x_p=\sum_{w\in W}d_w\min_{p\in P_w}\left\{\sum_{e\in E}\delta_{ep}\tilde{t}_e^0\right\}=\sum_{w\in W}d_w T_w\left(\tilde{t}^0\right),$$

*где* $T_w\left(\tilde{t}^0\right)$ *– длина кратчайшего пути из* $i$ *в* $j$ *(где* $w=(i,j)$*) на графе, ребра которого взвешены вектором* $\tilde{t}^0=\left\{\tilde{t}_e^0\right\}_{e\in E}$. *Таким образом, выписанную задачу можно решить с учетом того, что* $n=|V|\sim|E|$, *за* $\tilde{O}(Sn)$ *(здесь и далее* $\tilde{O}(\ )=O(\ )$ *с точностью до логарифмического фактора) и быстрее современными вариациями алгоритма Дейкстры [6, 23 – 25]. Обозначим решение этой задачи через* $f^0$.

Можно интерпретировать ситуацию таким образом, что в начальный момент водители посчитали, что все дороги абсолютно свободны и выбрали согласно этому предположению кратчайшие пути, соответствующие их целям, и поехали по этим маршрутам (путям). На практике более равномерное распределение водителей по путям в начальный момент может оказаться более предпочтительным.

Поняв, что в действительности из-за наличия других водителей время в пути не соответствует первоначальной оценке, выраженной весами ребер $\tilde{t}_e^0$, доля $\gamma^k$ водителей (обнаруживших это и готовых что-то менять) на следующем $(k+1)$-м шаге изменят свой выбор исходя из кратчайших путей, посчитанных по распределению водителей на предыдущем $k$-м шаге. Таким образом, возникает процедура “нащупывания” равновесия. Если выбирать специальным образом $\gamma^k$ (в частности, необходимо $\gamma^k\xrightarrow[k\to\infty]{}0$, чтобы избежать колебания вокруг равновесия (minority game [22]), и $\sum_{k=0}^{\infty}\gamma^k=\infty$, чтобы до равновесия дойти), то система, действительно, сойдется в равновесие. Опишем теперь более формально сказанное.

**Итерации** $k=0,1,2,\ldots$

*Пусть $f^k$ – вектор потоков на ребрах, полученный на предыдущей итерации с номером $k$. Положим $\tilde{t}_e^k = \partial\Psi\left(f^k\right)/\partial f_e = \tau_e\left(f^k\right)$ и рассмотрим задачу*

$$\sum_{e\in E}\tilde{t}_e^k y_e \to \min_{\substack{y=\Theta x\\ x\in X}} .$$

*Так же, как и раньше задача сводится к поиску кратчайших путей на графе, ребра которого взвешены вектором $\tilde{t}^k = \left\{\tilde{t}_e^k\right\}_{e\in E}$.*

*Обозначим решение задачи через $y^k$. Положим*

$$f^{k+1} = \left(1-\gamma^k\right)f^k + \gamma^k y^k,\ \gamma^k = \frac{2}{k+1}.$$

Заметим, что возникающую здесь задачу поиска кратчайших путей на графе можно попробовать (этому планируется посвятить отдельную работу) решать быстрее, чем за $\tilde{O}(Sn)$. Связано это с тем, что мы уже решали на предыдущей итерации аналогичную задачу для этого же графа с близкими весами ребер [6, 23 – 25] (веса ребер графа с ростом $k$ меняются все слабее от шага к шагу, поскольку $\gamma^k \xrightarrow[k\to\infty]{} 0$). Тем не менее, далее в статье мы будем считать, что одна итерация этого метода занимает $\tilde{O}(Sn)$.

Заметим также, что решая задачи поиска кратчайших путей мы находим (одновременно, т.е. без дополнительных затрат) не только вектор распределения потоков по ребрам $y$, но и разреженный вектор распределения потоков по путям $x$.

Строго говоря, нужно найти вектор $y^k$, а не кратчайшие пути. Чтобы получить вектор $y^k$ за $\tilde{O}(Sn)$ стоит для каждого из $S$ источников построить (например, алгоритмом Дейкстры) соответствующее дерево кратчайших путей (исходя из принципа динамического программирования “часть кратчайшего пути сама будет кратчайшим путем” несложно понять, что получится именно дерево, с корнем в рассматриваемом источнике). Это можно сделать для одного источника за $\tilde{O}(n)$. Однако, главное, правильно взвешивать ребра (их не больше $O(n)$) такого дерева, чтобы за один проход этого дерева можно было восстановить вклад (по всем ребрам) соответствующего источника в общий вектор $y^k$. Ребро должно иметь вес равный сумме всех проходящих через него корреспонденций с заданным источником (корнем дерева). Имея значения соответствующих корреспонденций (их также не больше $O(n)$) за один обратный проход (то есть с листьев к корню) такого дерева можно осуществить необходимое взвешивание (с затратами не более $O(n)$). Делается это по правилу: вес ребра равен сумме корреспонденции (возможно, равной нулю), в соответствующую вершину, в которую ребро входит и сумме весов всех ребер (если таковые имеются), выходящих из упомянутой вершины.

**Теорема 2 [18 – 21].** *Имеет место следующая оценка*

$$\Psi\left(f^N\right)-\Psi_* \le \Psi\left(f^N\right)-\Psi_N \le \frac{2L_p R_p^2}{N+1},\ f^N \in \Delta = \left\{ f = \Theta x:\ x \in X \right\},$$

*где*

$$\Psi_N = \max_{k=0,\dots,N}\left\{\Psi\left(f^k\right)+\left\langle \nabla\Psi\left(f^k\right), y^k - f^k\right\rangle\right\},$$

$$R_p^2 = \max_{k=0,\dots,N}\left\| y^k - f^k \right\|_p^2 \le \max_{f,\breve{f}\in\Delta}\left\|\breve{f} - f\right\|_p^2,\ L_p = \max_{\|h\|_p\le 1}\ \max_{f\in\mathrm{conv}\left(f^0, f^1,\dots,f^N\right)}\left\langle h, \mathrm{diag}\left\{\tau_e'\left(f_e\right)\right\} h\right\rangle,\ 1\le p\le\infty .$$

**Замечание 1.** Из доказательства этой теоремы [18 – 21] можно усмотреть немного более тонкий способ оценки $L_p$, в котором вместо $f \in \mathrm{conv}\left(f^0, f^1,\dots, f^N\right)$ можно брать

$$f \in \mathrm{conv}\left(f^0, f^1\right)\cup\mathrm{conv}\left(f^1, f^2\right)\cup\ldots\cup\mathrm{conv}\left(f^{N-1}, f^N\right).$$

Однако для небольшого упрощения выкладок мы будем использовать приведенный в формулировке теоремы огрубленный вариант.

К сожалению, в приложениях (см., пример о расщеплении потоков на личный и общественный транспорт в разделе 3) функции $\tau_e\left(f_e\right)$ могут иметь вертикальные асимптоты, что не позволяет равномерно по $N$ ограничить $L_p$ (даже если более тонко оценивать $L_p$, см. замечание 1). Такие случаи мы просто исключаем из рассмотрения для метода, описанного в этом разделе. Другими словами, мы считаем, что функции $\tau_e\left(f_e\right)$ заданы на положительной полуоси. К таким функциям относятся, например, BPR-функции (см. раздел 3).

Обратим внимание на то, что сам метод никак не зависит от выбора параметра $p$, от того какие получаются $R_p^2$ и $L_p$, в то время как оценка на число итераций, которые необходимо сделать для достижения заданной по функции (функционалу) точности, от этого выбора зависит. Как следствие, от этого выбора зависит и критерий останова (значение $\Psi_*$ нам априорно не известно.

Будем считать $p = 2$ (сопоставимые оценки, получаются и при выборе $p = \infty$):

$$L_2\left(f^0, f^1,\dots, f^N\right) = \max_{f\in\mathrm{conv}\left(f^0, f^1,\dots,f^N\right)}\ \max_{e\in E}\tau_e'\left(f_e\right) = \max_{e\in E}\tau_e'\left(\max_{k=0,\dots,N} f_e^k\right),\ R_2^2 = \max_{f,\breve{f}\in\Delta}\left\|\breve{f} - f\right\|_2^2 .$$

Величину $R_2^2$ мы можем оценить априорно (при это, к сожалению, получается довольно грубая оценка), т.е. можно считать её нам известной. Труднее обстоит дело с $L_2$. Далее предлагается оригинальный способ запуска метода Франк–Вульфа, критерий останова которого не требует априорного знания $L_2$.

Задаемся точностью $\varepsilon > 0$. Оцениваем $R_2^2$. Полагаем $L_2 = 1$ (для определенности). Запускаем метод Франк–Вульфа с $N(L_2) = 2L_2R_2^2/\varepsilon$. На каждом шаге проверяем условие (это делается за $\mathrm{O}(n)$)

$$L_2\left(f^0, f^1, ..., f^k\right) = \max_{e\in E} \tau'_e\left(\max_{l=0,...,k} f_e^l\right) \le L_2.$$

Если на всех шагах условие выполняется, то сделав $N(L_2)$ шагов, гарантированно получим решение с нужной точностью. Если же на каком-то шаге $k < N(L_2)$ условие нарушилось, т.е. $L_2\left(f^0, f^1, ..., f^k\right) > L_2$, то полагаем $L_2 := L_2\left(f^0, f^1, ..., f^k\right)$, пересчитываем $N(L_2)$ и переходим к следующему шагу. Таким образом, по ходу итерационного процесса мы корректируем критерий останова, оценивая необходимое число шагов по получаемой последовательности $\left\{f^k\right\}$. Специфика данной постановки, которая позволила так рассуждать, заключается в наличии явного представления

$$L_2\left(f^0, f^1, ..., f^k\right) = \max_{e\in E} \tau'_e\underbrace{\left(\max_{l=0,...,k} f_e^l\right)}_{\hat{f}_e},$$

и независимости используемого метода от выбора $L_2$ (шаг метода Франк–Вульфа $\gamma^k = 2(k+1)^{-1}$ вообще ни от каких параметров не зависит).

На практике, однако, приведенный способ работает не очень хорошо из-за использования завышенных оценок для $L_2$ и $R_2^2$. Более эффективным оказался другой способ, который использует неравенство (см. теорему 2) $\Psi\left(f^N\right) - \Psi_N \le \varepsilon$. В этом способе важно, что $\Psi_k$, $k = 0, ..., N$ – автоматически рассчитываются на каждой итерации без дополнительных затрат, а $\Psi\left(f^k\right)$ может быть рассчитано на каждой итерации по известному $f^k$ за $\tilde{\mathrm{O}}(n)$. Однако нет необходимости проверять этот критерий на каждой итерации, можно это делать, например, с периодом $\alpha\tilde{\varepsilon}^{-1}$, где $\tilde{\varepsilon}$ – относительная точность по функции (скажем, $\tilde{\varepsilon} = 0.01$ – означает, что $\varepsilon = 0.01\Psi\left(f^0\right)$), $\alpha \approx 1$ подбирается эвристически, исходя из задачи. Следуя [18], можно еще немного упростить рассуждения за счет небольшого увеличения числа итераций. А именно, можно использовать оценки (при этом следует полагать $\gamma^k = 2(k+2)^{-1}$)

$$\Psi\left(f^k\right) - \Psi_* \le \left\langle \nabla\Psi\left(f^k\right), f^k - y^k \right\rangle,$$

$$\min_{k=1,\dots,N} \left\langle \nabla\Psi\left(f^k\right), f^k - y^k \right\rangle \le \frac{7L_2R_2^2}{N+2}.$$

Таким образом, в данном разделе был описан способ поиска равновесного распределения потоков по ребрам $f$, который за время

$$\tilde{O}\left(SnL_2R_2^2/\varepsilon\right)$$

находит такой $f^{N(\varepsilon)}$, что

$$\Psi\left(f^{N(\varepsilon)}\right) - \Psi_* \le \varepsilon.$$

## 3. Рандомизированный метод двойственных усреднений поиска равновесия в модели стабильной динамики (Нестерова–де Пальма)

В ряде постановок задач вместо функций затрат на ребрах $\tau_e\left(f_e\right)$ заданы ограничения на пропускные способности $f_e \le \bar{f}_e$ и затраты на прохождения свободного (не загруженного $f_e < \bar{f}_e$) ребра $\bar{t}_e$. В модели стабильной динамики это сделано для всех ребер [1, 26], а в модели грузоперевозок РЖД – только для части [2]. Согласно работе [1], такую новую модель можно получить предельным переходом из модели Бэкмана, с помощью введения внутренних штрафов в саму модель. А именно, будем считать, что (как и в модели Бэкмана) у всех ребер есть свои функции затрат $\tau_e^\mu\left(f_e\right)$, но для части ребер $e \in E'$ (какой именно части, зависит от задачи) осуществляется предельный переход

$$\tau_e^\mu\left(f_e\right) \xrightarrow[\mu\to 0+]{} \begin{cases} \bar{t}_e, & 0 \le f_e < \bar{f}_e \\ \left[\bar{t}_e, \infty\right), & f_e = \bar{f}_e \end{cases},$$

$$d\tau_e^\mu\left(f_e\right)/df_e \xrightarrow[\mu\to 0+]{} 0,\ 0 \le f_e < \bar{f}_e.$$

Обозначив через $x(\mu)$ – равновесное распределение потоков по путям в модели Бэкмана при функциях затрат на ребрах $\tau_e^\mu\left(f_e\right)$, получим, что при $e \in E'$

$$\tau_e^\mu\left(f_e\left(x(\mu)\right)\right) \xrightarrow[\mu\to 0+]{} t_e,$$

$$f_e\left(x(\mu)\right) \xrightarrow[\mu\to 0+]{} f_e,$$

где пара $(t, f)$ – равновесие в модели стабильной динамики и ее вариациях [1, 2, 26] с тем же графом и матрицей корреспонденций, что и в модели Бэкмана, и с ребрами $e \in E'$, характеризующимися набором $(\bar{t}, \bar{f})$ из определения функций $\tau_e^{\mu}(f_e)$. Заметим, что если $t_e > \bar{t}_e$, то $t_e - \bar{t}_e$ можно интерпретировать, например, как время, потерянное в пробке на этом ребре [1, 26].

Согласно разделу 2 равновесная конфигурация при таком переходе $\mu \to 0+$ должна находиться из решения задачи

$$\Psi(f) = \sum_{e \in E \setminus E'} \int_0^{f_e} \tau_e^{\mu}(z)\,dz + \lim_{\mu \to 0+} \sum_{e \in E'} \int_0^{f_e} \tau_e^{\mu}(z)\,dz \to \min_{f = \Theta x,\, x \in X}.$$

Считая, что в равновесии не может быть $\lim_{\mu \to 0+} \tau_e^{\mu}(f_e) = \infty$ (иначе, равновесие просто не достижимо, и со временем весь граф превратится в одну большую пробку), можно не учитывать в интеграле вклад точек $\bar{f}_e$ (в случае попадания в промежуток интегрирования), то есть переписать задачу следующим образом

$$\min_{f = \Theta x,\, x \in X} \left\{ \sum_{e \in E \setminus E'} \int_0^{f_e} \tau_e^{\mu}(z)\,dz + \sum_{e \in E'} \int_0^{f_e} \left(\bar{t}_e + \delta_{\bar{f}_e}(z)\right) dz \right\} \Leftrightarrow \min_{\substack{f = \Theta x,\, x \in X \\ f_e \le \bar{f}_e,\, e \in E'}} \left\{ \sum_{e \in E \setminus E'} \int_0^{f_e} \tau_e^{\mu}(z)\,dz + \sum_{e \in E'} f_e \bar{t}_e \right\},$$

где

$$\delta_{\bar{f}_e}(z) = \begin{cases} 0, & 0 \le z < \bar{f}_e \\ \infty, & z \ge \bar{f}_e \end{cases}, \quad e \in E'.$$

**Теорема 3 [1, 26, раздел 4 ниже].** *Двойственная задача к выписанной выше задаче может быть приведена к следующему виду:*

$$\Upsilon(t) = -\sum_{w \in W} d_w T_w(t) + \left\langle \bar{f}, t - \bar{t} \right\rangle - \mu \sum_{e \in E \setminus E'} h_e^{\mu}(t_e) \to \min_{\substack{t_e \ge \bar{t}_e,\, e \in E' \\ t_e \in \operatorname{dom} h_e^{\mu}(t_e),\, e \in E \setminus E'}}, \qquad (1)$$

*где* $T_w(t)$ *– длина кратчайшего пути из* $i$ *в* $j$ *(*$w = (i, j) \in W$*) на графе, ребра которого взвешены вектором* $t = \{t_e\}_{e \in E}$*, а функции* $h_e^{\mu}(t_e)$ *– гладкие и вогнутые.*

*При этом решение изначальной задачи* $f$ *можно получить из формул:*

$$f_e = \bar{f}_e - s_e,\ e \in E',\ \text{где } s_e \ge 0 \text{ – множитель Лагранжа к ограничению } t_e \ge \bar{t}_e;$$

$$\tau_e^{\mu}(f_e) = t_e,\ e \in E \setminus E'.$$

Приведем пример модели (типа стабильной динамики) расщепления пользователей на личный и общественный транспорт [1], в которой каждое ребро $e \in E$ изначального графа продублировано для личного (“*л*”) и общественного (“*о*”) транспорта, при этом для общественного транспорта [26]

$$\tau_e\left(f_e^o\right)=\bar{t}_e^o\cdot\left(1+\mu\frac{\bar{f}_e^o}{\bar{f}_e^o-f_e^o}\right),$$

а для личного транспорта был осуществлен предельный переход $\mu\to 0+$ в аналогичных формулах

$$\tau_e\left(f_e^л\right)=\bar{t}_e^л\cdot\left(1+\mu\frac{\bar{f}_e^л}{\bar{f}_e^л-f_e^л}\right).$$

Поиск равновесного расщепления на [личный] и [общественный] транспорт приводит к следующей задаче [1]:

$$-\sum_{w\in W}d_w\min\left\{T_w^л\left(t^л\right),T_w^o\left(t^o\right)\right\}+\left\langle\bar{f}^л,t^л-\bar{t}^л\right\rangle+\left\langle\bar{f}^o,t^o-\bar{t}^o\right\rangle-\mu\sum_{e\in E}\bar{f}_e^o\cdot\bar{t}_e^o\cdot\ln\left(1+\frac{t_e^o-\bar{t}_e^o}{\bar{t}_e^o\mu}\right)\to\min_{\substack{t^л\ge\bar{t}^л\\ t^o\ge\bar{t}^o\cdot(1-\mu)}},$$

при этом $f^л=\bar{f}^л-s^л$, где $s^л$ – вектор множителей Лагранжа для ограничений $t^л\ge\bar{t}^л$,

$$f_e^o=\bar{f}_e^o\cdot\left(1-\frac{\bar{t}_e^o\cdot\mu}{t_e^o-(1-\mu)\bar{t}_e^o}\right).$$

Для упрощения рассуждений далее будем считать, что $E'=E$, т.е. на всех ребрах перешли к пределу $\mu\to 0+$. Если это не так, то надо будет далее использовать не метод зеркального спуска [19] (см. формулу (2) ниже), а его композитный вариант [27] с сепарабельным композитом $-\mu\sum_{e\in E\setminus E'}h_e^\mu\left(t_e\right)$. Возникающая на каждой итерации задача (поиска градиентного отображения, см. формулу (2)), в виду сепарабельности ограничений, распадается в $|E|=O(n)$ одномерных задач выпуклой оптимизации, которые можно решить за $\tilde{O}(n)$ с машинной точностью. Интересно заметить, что в случае, когда

$$\tau_e^\mu\left(f_e\right)=\bar{t}_e\cdot\left(1+\gamma\cdot\left(\frac{f_e}{\bar{f}_e}\right)^{\frac{1}{\mu}}\right)$$

– BPR-функции с $\mu=0.25$ (наиболее часто встречающаяся на практике ситуация [15]), то имеются явные формулы, поэтому итерацию можно сделать быстрее – за $O(n)$.

Численно решать задачу негладкой выпуклой оптимизации (1) с $E'=E$ будем с помощью специальным образом рандомизированных вариантов метода зеркального спуска [19, 28 – 31] с евклидовой прокс-структурой. Выбор такой прокс-структуры, прежде всего, связан с наличием ограничения $t\ge\bar{t}$ [28].

Пусть известно такое число $R$, что

$$R^2 \ge \frac{1}{2}\left\|t_* - t^0\right\|_2^2,$$

где $t_*$ – решение задачи (1). Выберем $N$ – число шагов алгоритма (далее, см. формулу (3), мы опишем, как можно выбирать $N$ исходя из желаемой точности $\varepsilon$ ). Положим начальное приближение $t^0 = \bar{t}$ . Пусть на шаге с номером $k$ получен вектор $t^k$ и стохастический субградиент $\nabla\Upsilon\left(t^k, \xi_k\right)$, зависящий от случайного вектора $\xi_k$, и удовлетворяющий условию несмещенности $E_{\xi_k}\left[\nabla\Upsilon\left(t^k, \xi_k\right)\right] = \nabla\Upsilon\left(t^k\right)$. Здесь и далее под $\nabla\Upsilon(t)$ мы имеем в виду какой-то измеримый селектор [32] многозначного отображения $\partial\Upsilon(t)$ (субдифференциала функции $\Upsilon(t)$). Следующая точка вычисляется из соотношения

$$t^{k+1} = \arg\min_{t \ge \bar{t}} \Big\{ \underbrace{\Upsilon\left(t^k\right)}_{\text{можно не писать}} + \underbrace{\gamma_k \left\langle \nabla\Upsilon\left(t^k, \xi_k\right), t - t^k \right\rangle}_{\nabla\Upsilon\left(t^k, \xi_k\right) \text{ - стохастический субградиент}} + \frac{1}{2}\left\|t - t^k\right\|_2^2 \Big\}, \quad k = 0, \ldots, N. \qquad (2)$$

Здесь

$$\gamma_k = \frac{R}{M_k}\sqrt{\frac{2}{N+1}},$$

где

$$\max\left\{\left\|\nabla\Upsilon\left(t^k, \xi_k\right)\right\|_2, \left\|\nabla\Upsilon\left(t^k\right)\right\|_2\right\} \le M_k \le M.$$

Отметим, что такой выбор $\gamma_k$ обусловлен решением задачи

$$2R^2/\tilde{\gamma} + (N+1)M^2\tilde{\gamma} \to \min_{\tilde{\gamma} \ge 0}, \text{ где } M = \max_{k=0,\ldots,N} M_k,$$

возникающей при доказательстве теоремы 4, см. ниже.

В виду сепарабельности ограничений $t \ge \bar{t}$ и сепарабельности выражения, стоящего в фигурных скобках в (2), задача (2) на каждом шаге итерационного процесса декомпозируется на $|E| = O(n)$ одномерных подзадач, каждая из которых представляет собой задачу минимизации параболы на полуоси. Следовательно, каждая такая подзадача решается по явным формулам, т.е. на каждой итерации за $O(n)$ можно решить задачу (2) в предположении, что мы нашли стохастический субградиент $\nabla\Upsilon\left(t^k, \xi_k\right)$. Чтобы оценить насколько много потребуется итераций $N = N(\varepsilon)$ для достижения заданной точности $\Upsilon\left(\bar{t}^{N(\varepsilon)}\right) - \Upsilon_* \le \varepsilon$, где $\Upsilon_* = \min_{t \ge \bar{t}} \Upsilon(t)$, сформулируем теорему о сходимости метода (2).

**Теорема 4 [19, 28 – 31].** *Пусть*

$$\bar{t}^N = \frac{1}{S_N}\sum_{k=0}^{N}\gamma_k t^k, \quad S_N = \sum_{k=0}^{N}\gamma_k, \quad \bar{f}^N = \bar{f} - s^N,$$

*где* $s^N$ *– есть множитель Лагранжа к ограничению* $t \ge \bar{t}$ *в задаче*

$$\frac{1}{S_N}\left\{\sum_{k=0}^{N}\gamma_k\left\{\Upsilon\left(t^k\right)+\left\langle\nabla\Upsilon\left(t^k,\xi_k\right),t-t^k\right\rangle\right\}+\frac{1}{2}\left\|t-t^0\right\|_2^2\right\}\to\min_{t\ge\bar{t}}.$$

*Тогда с вероятностью* $\ge 1-\sigma$

$$0 \le \Upsilon\left(\bar{t}^N\right)-\Upsilon_* \le \frac{16\sqrt{2}MR}{\sqrt{N}}\ln\left(\frac{4N}{\sigma}\right)$$

*с*

$$R^2 = \frac{1}{2}\left\|t_* - t^0\right\|_2^2$$

*и с вероятностью* $\ge 1-\sigma$

$$0 \le \Upsilon\left(\bar{t}^N\right)+\Psi\left(\bar{f}^N\right) \le \frac{16\sqrt{2}MR}{\sqrt{N}}\ln\left(\frac{4N}{\sigma}\right) \qquad (3)$$

*с*

$$R^2 \ge \max\left\{\frac{1}{2}\left\|t_* - t^0\right\|_2^2, \frac{1}{2}\left\|\tilde{t}^N - t^0\right\|_2^2\right\},$$

$$\tilde{t}^N = \arg\max_{t\ge 0}\left\{\sum_{w\in W} d_w T_w(t) - \left\langle \bar{f}^N, t-\bar{t}\right\rangle\right\}.$$

**Доказательство.** Согласно [19, п. 3.4 29]

$$0 \le \Upsilon\left(\bar{t}^N\right)-\Upsilon_* \le$$

$$\le \frac{1}{S_N}\left\{\frac{1}{2}\left\|t_*-t^0\right\|_2^2-\frac{1}{2}\left\|t_*-t^{N+1}\right\|_2^2-\sum_{k=0}^{N}\gamma_k\left\langle\nabla\Upsilon\left(t^k,\xi_k\right)-\nabla\Upsilon\left(t^k\right),t^k-t_*\right\rangle+\frac{1}{2}\sum_{k=0}^{N}\gamma_k^2 M_k^2\right\}\le$$

$$\le \frac{R^2}{S_N}-\frac{1}{2S_N}\left\|t_*-t^{N+1}\right\|_2^2+\frac{1}{S_N}\sum_{k=0}^{N}\gamma_k\left\langle\nabla\Upsilon\left(t^k,\xi_k\right)-\nabla\Upsilon\left(t^k\right),t_*-t^k\right\rangle+\frac{R^2}{S_N}\le$$

$$\le \frac{\sqrt{2}MR}{\sqrt{N+1}}+\frac{1}{S_N}\sum_{k=0}^{N}\gamma_k\left\langle\nabla\Upsilon\left(t^k,\xi_k\right)-\nabla\Upsilon\left(t^k\right),t_*-t^k\right\rangle. \qquad (4)$$

При формировании последнего слагаемого предпоследнего неравенства было учтено определение $\gamma_k$.

Считая, что $\frac{1}{2}\left\|t_*-t^k\right\|_2^2 \le R_\sigma^2$ для всех $k=0,...,N$ с вероятностью $\ge 1-\sigma/2$, получим из неравенства Азума–Хефдинга [33, 34] для ограниченной мартингал–разности

$$\left|\left\langle \nabla\Upsilon\left(t^k,\xi_k\right)-\nabla\Upsilon\left(t^k\right),t_*-t^k\right\rangle\right|\le\left\|\nabla\Upsilon\left(t^k,\xi_k\right)-\nabla\Upsilon\left(t^k\right)\right\|_2\left\|t_*-t^k\right\|_2\le 2\sqrt{2}M_kR_\sigma$$

следующее неравенство

$$P\left(\sum_{k=0}^N\gamma_k\left\langle \nabla\Upsilon\left(t^k,\xi_k\right)-\nabla\Upsilon\left(t^k\right),t_*-t^k\right\rangle\ge 2\sqrt{2}R_\sigma\Lambda\sqrt{\sum_{k=0}^N\gamma_k^2M_k^2}\right)\le\exp\left(-\Lambda^2/2\right)+\sigma/2.$$

Следовательно,

$$P\left(\sum_{k=0}^N\gamma_k\left\langle \nabla\Upsilon\left(t^k,\xi_k\right)-\nabla\Upsilon\left(t^k\right),t_*-t^k\right\rangle\ge 4\sqrt{2}R_\sigma R\sqrt{\ln\left(2/\sigma\right)}\right)\le\sigma,$$

$$P\left(\frac{1}{S_N}\sum_{k=0}^N\gamma_k\left\langle \nabla\Upsilon\left(t^k,\xi_k\right)-\nabla\Upsilon\left(t^k\right),t_*-t^k\right\rangle\ge \frac{4MR_\sigma\sqrt{\ln\left(2/\sigma\right)}}{\sqrt{N+1}}\right)\le\sigma. \qquad (5)$$

Из неравенства (4) имеем

$$\frac{1}{2}\left\|t_*-t^k\right\|_2^2\le 2R^2+\sum_{l=0}^{k-1}\gamma_l\left\langle \nabla\Upsilon\left(t^l,\xi_l\right)-\nabla\Upsilon\left(t^l\right),t_*-t^l\right\rangle.$$

Неравенство (5) представим в виде $k=1,...,N$

$$P\left(\sum_{l=0}^{k-1}\gamma_l\left\langle \nabla\Upsilon\left(t^l,\xi_l\right)-\nabla\Upsilon\left(t^l\right),t_*-t^l\right\rangle\ge 4\sqrt{2}R_\sigma R\sqrt{\ln\left(2/\sigma\right)}\right)\le\sigma.$$

Отсюда по неравенству Буля (вероятность суммы событий, соответствующих $k=1,...,N$, не больше суммы вероятностей событий) имеем с вероятностью $\ge 1-\sigma$ для всех $k=0,...,N$

$$\frac{1}{2}\left\|t_*-t^k\right\|_2^2\le 2R^2+4\sqrt{2}R_\sigma R\sqrt{\ln\left(2N/\sigma\right)}.$$

Положим $R_\sigma^2=2R^2+4\sqrt{2}R_\sigma R\sqrt{\ln\left(4N/\sigma\right)}$. Отсюда получаем следующую оценку

$$R_\sigma\le 4R\left(1+\sqrt{2\ln\left(4N/\sigma\right)}\right).$$

Следовательно, из неравенств (4), (5) имеем с вероятностью $\ge 1-\sigma$ неравенство

$$0\le\Upsilon\left(\bar{t}^N\right)-\Upsilon_*\le\frac{\sqrt{2}MR}{\sqrt{N+1}}\left(1+8\sqrt{2}\left(1+\sqrt{2\ln\left(4N/\sigma\right)}\right)\sqrt{\ln\left(2/\sigma\right)}\right)\le\frac{16\sqrt{2}MR}{\sqrt{N}}\ln\left(\frac{4N}{\sigma}\right).$$

Положим теперь (см. также конец п. 3 статьи [35])

$$R^2:=\max\left\{\frac{1}{2}\left\|t_*-t^0\right\|_2^2,\frac{1}{2}\left\|\tilde{t}^N-t^0\right\|_2^2\right\},$$

где

$$\tilde{t}^N=\arg\max_{t\ge 0}\left\{\sum_{w\in W}d_wT_w\left(t\right)-\left\langle\bar{f}^N,t-\bar{t}\right\rangle\right\}.$$

Тогда

$$\Upsilon\left(\bar{t}^N\right) \le \frac{1}{S_N}\sum_{k=0}^{N}\gamma_k \Upsilon\left(t^k\right) \le \frac{1}{2S_N}\sum_{k=0}^{N}\gamma_k^2 M_k^2 +$$

$$+\frac{1}{S_N}\min_{t\ge\bar{t}}\left\{\sum_{k=0}^{N}\gamma_k\left\{\Upsilon\left(t^k\right)+\left\langle\nabla\Upsilon\left(t^k,\xi_k\right),t-t^k\right\rangle\right\}+\frac{1}{2}\left\|t-t^0\right\|_2^2\right\} \le \frac{1}{2S_N}\sum_{k=0}^{N}\gamma_k^2 M_k^2 +$$

$$+\min_{t}\left\{\frac{1}{S_N}\left\{\sum_{k=0}^{N}\gamma_k\left\{\Upsilon\left(t^k\right)+\left\langle\nabla\Upsilon\left(t^k,\xi_k\right),t-t^k\right\rangle\right\}+\frac{1}{2}\left\|t-t^0\right\|_2^2\right\}+\left\langle s^N,\bar{t}-t\right\rangle\right\} \le$$

$$\le \frac{1}{2S_N}\sum_{k=0}^{N}\gamma_k^2 M_k^2 + \min_{t\ge 0}\left\{\frac{1}{S_N}\left\{\sum_{k=0}^{N}\gamma_k\left\langle\nabla\Upsilon\left(t^k,\xi_k\right)-\nabla\Upsilon\left(t^k\right),t-t^k\right\rangle+\frac{1}{2}\left\|t-t^0\right\|_2^2\right\}+\right.$$

$$\left.+\Upsilon(t)+\left\langle s^N,\bar{t}-t\right\rangle\right\} \le$$

$$\le \frac{1}{2S_N}\sum_{k=0}^{N}\gamma_k^2 M_k^2 + \frac{1}{S_N}\sum_{k=0}^{N}\gamma_k\left\langle\nabla\Upsilon\left(t^k,\xi_k\right)-\nabla\Upsilon\left(t^k\right),\tilde{t}^N-t^k\right\rangle + \frac{1}{2S_N}\left\|\tilde{t}^N-t^0\right\|_2^2 - \Psi\left(\bar{f}^N\right).$$

К сожалению, в данном случае неравенство Азума–Хефдинга уже не будет выполняться, поскольку $\tilde{t}^N$ зависит от $\xi_k$ – нарушается условие, что последовательность $\left\{\gamma_k\left\langle\nabla\Upsilon\left(t^k,\xi_k\right)-\nabla\Upsilon\left(t^k\right),\tilde{t}^N-t^k\right\rangle\right\}$ будет последовательностью мартингал-разностей. Тем не менее, аналогично рассуждая, можно показать, что существует такое

$$R^2 \ge \max\left\{\frac{1}{2}\left\|t_*-t^0\right\|_2^2, \frac{1}{2}\left\|\tilde{t}^N-t^0\right\|_2^2\right\},$$

что с вероятностью $\ge 1-\sigma$

$$0 \le \Upsilon\left(\bar{t}^N\right)+\Psi\left(\bar{f}^N\right) \le \frac{16\sqrt{2}MR}{\sqrt{N}}\ln\left(\frac{4N}{\sigma}\right). \quad \square$$

**Замечание 2.** Выписанная в теореме 4 оценка (3) была получена с помощью теоремы 3 работы [29] и замечания 4 работы [36]. Более аккуратный способ рассуждений позволяет немного улучшить оценку (3). Для практических приложений метода (2) удобнее записать его не в терминах неизвестного $R$, а сразу в терминах другого неизвестного

$$\bar{R}^2 \ge \max\left\{\max_{k=0,\ldots,N}\left\|t_*-t^k\right\|_2^2, \frac{1}{2}\left\|\tilde{t}^N-t^0\right\|_2^2\right\}.$$

Речь идет о выборе шагов

$$\gamma_k = \frac{\bar{R}}{M_k}\sqrt{\frac{2}{N+1}},$$

и об оценке [31]

$$0 \le \Upsilon\left(\bar{t}^{N}\right)+\Psi\left(\bar{f}^{N}\right) \le \frac{\sqrt{2}M\bar{R}}{\sqrt{N}}\left(1+\sqrt{8\ln\left(2/\sigma\right)}\right)$$

с вероятностью $\ge 1-\sigma$. С практической точки зрения от этого ничего не меняется, и даже становится лучше (оценка улучшается). Этот параметр как был априорно неизвестным, так и остался таковым, поменялась немного только его интерпретация. Тем не менее, оценка (3) представляет определенный теоретический интерес, поскольку является хорошей демонстрацией тех теоретических трудностей, которые возникают при описанном выше способе восстановления решения прямой задачи.

**Замечание 3.** Пусть $\nabla\Upsilon\left(t^{k},\xi_{k}\right)\equiv\nabla\Upsilon\left(t^{k}\right)$, т.е. вместо стохастического субградиента в алгоритме (2) используется обычный субградиент. Тогда вектор равновесного распределения потоков по ребрам $f$ можно считать по-другому (при таком способе подсчета возможно нарушение неравенств $\bar{f}^{N} \le \bar{f}$ и $0 \le \Upsilon\left(\bar{t}^{N}\right)+\Psi\left(\bar{f}^{N}\right)$):

$$f^{k} \in \partial \sum_{w\in W} d_{w}T_{w}\left(t^{k}\right),\ \bar{f}^{N} = \frac{1}{S_{N}}\sum_{k=0}^{N}\gamma_{k} f^{k}.$$

Решая возникающие здесь задачи поиска кратчайших путей, мы находим (одновременно, т.е. без дополнительных затрат) не только вектор распределения потоков по ребрам $f^{k}$, но и разреженный вектор распределения потоков по путям $x^{k}$ (при этом $f^{k} = \Theta x^{k}$). При этом

$$\Upsilon\left(\bar{t}^{N}\right)+\Psi\left(\bar{f}^{N}\right)+2\sqrt{2}R\left\|\left(\bar{f}^{N}-\bar{f}\right)_{+}\right\|_{2} \le \frac{2\sqrt{2}MR}{\sqrt{N}},$$

где

$$R^{2} = \frac{1}{2}\left\|t_{*}-t^{0}\right\|_{2}^{2}.$$

Отсюда можно получить, что

$$0 \le \Upsilon\left(\bar{t}^{N}\right)-\Upsilon_{*} \le \frac{\sqrt{2}MR}{\sqrt{N}},$$

$$\left|\Psi\left(\bar{f}^{N}\right)-\Psi_{*}\right| \le \frac{2\sqrt{2}MR}{\sqrt{N}},$$

$$\left\|\left(\bar{f}^{N}-\bar{f}\right)_{+}\right\|_{2} \le \frac{2M}{\sqrt{N}}.$$

Из теоремы 4 и слабой двойственности ($-\Upsilon_{*} \le \Psi_{*}$) имеем

$$\Upsilon\left(\bar{t}^{N}\right)-\Upsilon_{*}+\Psi\left(\bar{f}^{N}\right)-\Psi_{*} \le \Upsilon\left(\bar{t}^{N}\right)+\Psi\left(\bar{f}^{N}\right) \le \varepsilon.$$

Это обосновывает следующее следствие теоремы 4.

**Следствие.** *В условиях теоремы 4*

$$0 \le \Upsilon\left(\bar{t}^N\right) - \Upsilon_* \le \varepsilon,\ 0 \le \Psi\left(\bar{f}^N\right) - \Psi_* \le \varepsilon.$$

Осталось описать, как можно случайно выбирать быстро вычислимый, равномерно ограниченный по норме, стохастический субградиент функции $\Upsilon(t)$ со свойством несмещенности.

Прежде всего, опишем, как формируется субградиент функции $\Upsilon(t)$. Заметим, что субградиент выпуклой негладкой функции $\partial\Upsilon(t)$ – есть выпуклое множество, превращающееся в точках гладкости $\Upsilon(t)$ в один вектор – обычный градиент $\nabla\Upsilon(t)$. Из определения $\Upsilon(t)$ (см. формулу (1)) имеем

$$\partial\Upsilon(t) = -\sum_{w\in W} d_w \partial T_w(t) + \bar{f},$$

где $\partial T_w(t)$ – супердифференциал негладкой вогнутой (как минимум выпуклых, в нашем случае аффинных) функции $T_w(t) = \min_{p\in P_w} \sum_{e\in E} \delta_{ep} t_e$ (следует сравнить с задачами, возникающими на каждом шаге метода Франк–Вульфа из раздела 2). Супердифференциал $\partial T_w(t)$ представляет собой выпуклую комбинацию векторов с компонентами $\delta_{ep}$, отвечающих кратчайшим путям (если их несколько) для заданной корреспонденции $w$ на графе, ребра которого взвешены вектором $t = \{t_e\}_{e\in E}$. Каждый такой вектор (с числом компонент равным числу ребер) можно описать следующем образом: если ребро входит в кратчайший путь, то в компоненте вектора, отвечающего этому ребру, стоит 1, иначе 0.

Теперь опишем два варианта (первый вариант был сообщен нам Ю.Е. Нестеровым в 2013 г.) выбора несмещенного стохастического субградиента (мы вводим случайность, говорят также рандомизацию, чтобы за счет этого сократить стоимость вычисления)

**Вариант 1**

$$\nabla\Upsilon(t,\xi) = -d\nabla T_\xi(t) + \bar{f},$$

*где с.в.* $\xi = w$ *с вероятностью* $d_w/d$, $w \in W$, $\nabla T_\xi(t)$ *– произвольный элемент супердифференциала* $\partial T_\xi(t)$.

**Вариант 2**

$$\nabla\Upsilon(t,\xi) = -d \sum_{j:(\xi,j)\in W} \frac{d_{\xi j}}{d_{\xi\cdot}} \nabla T_{\xi j}(t) + \bar{f},$$

*где с.в.* $\xi = i$ *с вероятностью* $d_{i\cdot}/d$, $d_{i\cdot} = \sum_{j:(i,j)\in W} d_{ij}$, $i \in V$, $\nabla T_{\xi}(t)$ *– произвольный элемент супердифференциала* $\partial T_{\xi}(t)$.

Заметим, что в варианте 2 $\nabla\Upsilon(t,\xi)$ может быть вычислен алгоритмом Дейкстры за $\tilde{O}(n)$ в силу особенность алгоритма Дейкстры поиска кратчайших путей, заключающейся в том, что за время $\tilde{O}(n)$ он находит кратчайшие пути из заданной вершины во все остальные [6]. Однако подсчет суммы в определении $\nabla\Upsilon(t,\xi)$, если это делать напрямую, может занять время $O(n^{3/2})$, поскольку число слагаемых равно $O(n)$, а число ненулевых компонент в векторах $\nabla T_{\xi j}(t)$ (число ребер в соответствующем кратчайшем пути) может быть порядка $O(\sqrt{n})$ (сеть типа двумерной решетки). Однако можно по-другому организовать вычисление компонент вектора $\nabla\Upsilon(t,\xi)$. Алгоритм Дейкстры строит (ориентированное) дерево кратчайших путей с корнем в $\xi$ за $\tilde{O}(n)$. Для каждой вершины этого дерева $j \neq \xi$, припишем ребру, ведущему в эту вершину вес $d_{\xi j}$. Получим таким образом взвешенное дерево. Далее припишем листьям дерева веса, равные весам ребер, ведущим в эти листья. А далее по индукции: припишем вершине дерева сумму весов всех потомков. Припишем теперь ребрам дерева новые веса: вес ребра равен весу вершины, в которую это ребро входит. Далее нужно просто пробежаться по ребрам этого дерева, считывая веса и нормируя их на $-d/d_{\xi\cdot}$. С точностью до $\bar{f}$ получим таким образом $\nabla\Upsilon(t,\xi)$. Все это можно сделать за $\tilde{O}(n)$. Аналогичные конструкции позволяют вычислять $\nabla\Upsilon(t)$ и $\Upsilon(\bar{t}^N)$ за $\tilde{O}(Sn)$.

Заметим также, что в варианте 2 оценка константы $M$ (см. теорему 4) получается заметно лучше, чем в варианте 1.

Таким образом, вариант 2 кажется более предпочтительным. В чем может быть минус использования варианта 2 – в возможности использовать более быстрые алгоритмы поиска кратчайших путей, которые находят кратчайший путь ровно между двумя вершинами. Хотя нижняя оценка затрат на поиск кратчайшего пути здесь также $\tilde{O}(n)$ для реальных (почти планарных) транспортных (и не только) сетей эта оценка может быть существенно редуцирована. В частности, в определенных ситуациях до $\tilde{O}(1)$. Как правило, это, в свою очередь, требует затрат порядка $\tilde{O}(n)$ на подготовку специальной структуры данных [6] (будем называть такой процесс препроцессингом). Тем не менее, учитывая, что веса ребер меняются от шага к шагу не сильно, такой препроцессинг не обязательно осуществлять на каждой итерации. К тому же совершенно не обязательно находить всегда кратчайшие пути, считая таким образом точный субдифференциал. Желая решить задачу с точностью $\varepsilon$ по зазору двойственности $\Upsilon\left(\bar{t}^N\right)+\Psi\left(\bar{f}^N\right)\le\varepsilon$ достаточно вычислять $O(\varepsilon)$-субдифференциал на каждой итерации [37].

Выше мы пояснили, что нельзя однозначно, исходя из теоретических оценок, отдать предпочтение одному из описанных вариантов выбора $\nabla\Upsilon(t,\xi)$. Более того, в контексте проводимого выше анализа методов может быть не очевидна и сама необходимость в рандомизации.

Действительно, на первый взгляд кажется, что если считать полностью $\nabla\Upsilon(t)$ (такой подход с методом двойственных усреднений [29] вместо зеркального спуска предлагался Ю.Е. Нестеровым в 2012 г.), то оценка (3) примет вид (заметим, что с точностью до множителя $\sqrt{2}$ эта оценка – не улучшаема [19, 28, 37])

$$\Upsilon\left(\bar{t}^N\right)+\Psi\left(\bar{f}^N\right)\le\frac{\sqrt{2}MR}{\sqrt{N}},$$

Поскольку вычисление $\nabla\Upsilon(t)$ требует $\tilde{O}(Sn)$ операций (в определенных ситуациях возможно и быстрее), то, кажется, что при таком подходе мы просто теряем фактор $n$ в оценке сложности метода. Однако в действительности это не совсем так. Во-первых, константа $M$ здесь может быть заметно меньше своего аналога в варианте 2. Во-вторых, поскольку на каждой итерации мы должны пересчитывать все кратчайшие пути, то в этой постановке также как и в разделе 2 на передний план выходит пересчет кратчайших путей вместо расчета, что с учетом допустимости использования приближенно вычисленных кратчайших путей может редуцировать оценку сложности итерации $\tilde{O}(Sn)$ [6, 23 – 25].

Далее мы, тем не менее, ограничимся рассмотрением рандомизированных методов, считая стоимость одной итерации равной $\tilde{O}(n)$.

Свойство несмещенности стохастических субградиентов следует из построения (вариант 1 и 2). Сложнее обстоит дело с определением параметра метода $R$ (или $R_\sigma$, см. замечание 2; далее мы ограничимся рассмотрением случая, когда в качестве параметра выбран $R$, с помощью замечания 2 можно провести аналогичные рассуждения и для $R_\sigma$), который явно входит в итерационный процесс (2) (следует сравнить с методом Франк–Вульфа, для которого параметры метода не входили в сам метод, только в один из вариантов критерия останова).

Критерий останова можно задавать явной формулой для числа итераций (см. теорему 4 и замечание 2), в которую входит неизвестный параметр $R$ (или $\bar{R}$), но лучше его задавать немного по-другому (см. ниже).

Прежде всего, отметим, что если известна оценка (сверху) на $M$, то ее можно использовать при выборе шага метода

$$\gamma_k \equiv \frac{R}{M}\sqrt{\frac{2}{N+1}} \text{ или } \gamma_k \equiv \frac{\bar{R}}{M}\sqrt{\frac{2}{N+1}},$$

при этом теорема 4 и замечание 2 останутся верными. Далее сконцентрируемся на замечании 2. В формуле для $\gamma_k$ стоит неизвестное $\bar{R}$, которое исчезает при подстановке сюда зависимости $N\left(\varepsilon;\bar{R},M,\sigma\right)$, определяемой в замечании 2 из условия

$$\frac{\sqrt{2}M\bar{R}}{\sqrt{N\left(\varepsilon;\bar{R},M,\sigma\right)}}\left(1+\sqrt{8\ln\left(2/\sigma\right)}\right)=\varepsilon .$$

Таким образом, шаг метода не зависит от неизвестного $\bar{R}$. В качестве, критерия остановки метода используется условие (проверяемое за $\tilde{O}\left(Sn\right)$)

$$\Upsilon\left(\bar{t}^N\right)+\Psi\left(\bar{f}^N\right)\le\varepsilon .$$

Если априорно оценка константы $M$ не известна, то процедура усложняется. Сначала предполагается, что $M=M_0$. Далее предлагается следующая процедура адаптивного подбора неизвестного параметра $\bar{R}$ (отчасти являющаяся оригинальной, см., например, [36]). Задаем какое-то начальное значение, скажем, $\bar{R}=\left\|t^0\right\|_2$, запускаем итерационный процесс (2) с этими параметрами. В какой-то момент мы либо обнаружим, что $M>M_0$, либо будет сделано предписанное (теоремой 4) для выбранных значений параметров число итераций. Предположим, что имеет место вторая альтернатива. Далее проверяем условие

$$\Upsilon\left(\bar{t}^N\right)+\Psi\left(\bar{f}^N\right)\le\varepsilon .$$

Если оно выполняется, то мы нашли решение с требуемой точностью. Если не выполняется, то запускаем процесс заново, полагая $\bar{R}:=\sqrt{2}\bar{R}$ (такой выбор константы также оптимален, см. [38]). Это дополнительно может привести к не более чем логарифмическому (от отношения истинного значения $\bar{R}$ к $\left\|t^0\right\|_2$) числу перезапусков. Здесь в рассуж-дениях мы пренебрегли вероятностными оговорками, поскольку $\sigma$ можно считать очень малым (см. формулу (3)).

Если мы вышли из описанного цикла из-за того, что на каком-то шаге получили $M > M_0$, то полагаем $M := \sqrt{2}M$ (такой выбор константы также оптимален, см. [38]) и запускаем итерационный процесс заново с новым значением $M$. Это дополнительно может привести к не более чем логарифмическому (от отношения истинного значения $M$ к $M_0$) числу перезапусков.

В итоге, ожидаемая оценка времени работы метода (2) – есть $\tilde{O}\left(nM^2R^2/\varepsilon^2\right)$.

## 4. Рандомизированный покомпонентный спуск для модели стабильной динамики в форме Ю.Е. Нестерова

Предположим, что число источников намного меньше числа вершин

$$S = |O| \ll |V| = n .$$

Описанные в разделах 2, 3 методы не сильно учитывают такую разреженность. В мае 2014 года Ю.В. Дорном и Ю.Е. Нестеровым была предложена следующая эквивалентная переформулировка (см. Теорему 5 в Приложение) двойственной задачи (1) для модели стабильной динамики [26, 39] (т.е. при $E' = E$)

$$\min_{t \ge \bar{t}} \left\{ -\sum_{w \in W} d_w T_w(t) + \left\langle \bar{f}, t - \bar{t} \right\rangle \right\} = \min_{T} \left\{ -\sum_{\substack{s \in O,\, k \in D \\ (s,k) \in W}} d_{sk} \cdot \left(T_{sk} - T_{ss}\right) + \sum_{(i,j) \in E} \bar{f}_{ij} \max_{s \in O} \left(T_{sj} - T_{si} - \bar{t}_{ij}, 0\right) \right\}, \quad (6)$$

которая, по-видимому, позволяет в большей степени учесть свойство $S \ll n$. По решению задачи (6) можно восстанавливать решение двойственной задачи (1):

$$t_{ij} = \max\left\{ \max_{s \in O}\left(T_{sj} - T_{si}\right), \bar{t}_{ij} \right\},\ (i,j) \in E ,$$

но, к сожалению, нельзя восстанавливать вектор равновесного распределения потоков по путям $x$. Учитывая что

$$\min_{f \in \Delta,\, f \le \bar{f}} \sum_{e \in E} f_e \bar{t}_e = \min_{\substack{f = \Theta x,\, x \in X \\ f \le \bar{f}}} \sum_{e \in E} f_e \bar{t}_e = \min_{f = \Theta x,\, x \in X} \max_{\tau \ge 0} \left\{ \sum_{e \in E} f_e \cdot \left(\bar{t}_e + \tau_e\right) - \left\langle \bar{f}, \tau \right\rangle \right\} =$$

$$= \max_{\tau \ge 0} \left\{ \sum_{w \in W} d_w T_w\left(\bar{t} + \tau\right) - \left\langle \bar{f}, \tau \right\rangle \right\} \overset{t = \bar{t} + \tau}{=} -\min_{t \ge \bar{t}} \left\{ -\sum_{w \in W} d_w T_w(t) + \left\langle \bar{f}, t - \bar{t} \right\rangle \right\},$$

можно получить по формуле Демьянова–Данскина–Рубинова [1, 28] и решение прямой задачи

$$f\in\partial_{\bar t}\left(\min_{f\in\Delta,\, f\le\bar f}\sum_{e\in E} f_e\bar t_e\right)=\partial_{\bar t}\left(-\min_T\left\{-\sum_{\substack{s\in O,\,k\in D\\(s,k)\in W}} d_{sk}\cdot\left(T_{sk}-T_{ss}\right)+\sum_{(i,j)\in E}\bar f_{ij}\max_{s\in O}\left(T_{sj}-T_{si}-\bar t_{ij},0\right)\right\}\right).$$

Следуя [39, 40], с помощью техники двойственного сглаживания, запишем функционал, равномерно аппроксимирующий целевой функционал задачи (6) (здесь $\eta_{ij}=\varepsilon/\left(4n\bar f_{ij}\ln\left(S+1\right)\right)$, $\varepsilon$ – точность, с которой хотим решить задачу (1)) в виде

$$-\sum_{s\in O,k\in D} d_{sk}\cdot\left(T_{sk}-T_{ss}\right)+\sum_{(i,j)\in E}\bar f_{ij}\eta_{ij}\ln\left(\frac{1}{S+1}\left[\sum_{s\in O}\exp\left(\frac{T_{sj}-T_{si}-\bar t_{ij}}{\eta_{ij}}\right)+1\right]\right)\to\min_T, \qquad (7)$$

$$f_{ij}=\bar f_{ij}\frac{\sum_{s\in O}\exp\left(\left(T_{sj}-T_{si}-\bar t_{ij}\right)/\eta_{ij}\right)}{\sum_{s\in O}\exp\left(\left(T_{sj}-T_{si}-\bar t_{ij}\right)/\eta_{ij}\right)+1}=\frac{\bar f_{ij}}{1+\left(\sum_{s\in O}\exp\left(\left(T_{sj}-T_{si}-\bar t_{ij}\right)/\eta_{ij}\right)\right)^{-1}},\ (i,j)\in E.$$

Решая задачу (7) с точностью $\varepsilon/4$ по функции прямо-двойственным методом (недавно было установлено [41], что при правильном взгляде любой разумный численный метод является прямо-двойственным, точнее имеет соответствующую модификацию) можно восстановить $\tilde t^N$ и $\tilde f^N$ так, чтобы зазор двойственности был меньше $\varepsilon$

$$\Upsilon\left(\tilde t^N\right)+\Psi\left(\tilde f^N\right)\le\varepsilon.$$

Таким образом, достаточно научиться эффективно решать задачу (7). Это можно сделать, например, с помощью ускоренного покомпонентного спуска Ю.Е. Нестерова, или более современных вариантов ускоренных покомпонентных спусков APPROX, ALPHA [42 – 44]. Разобьем все компоненты $\left\{T_{sk}\right\}_{s\in O,k\in V}$ на блоки $\left\{T_{(k)}\right\}_{k\in V}$, где $T_{(k)}=\left\{T_{sk}\right\}_{s\in O}$. Чтобы получить оценку скорости сходимости, нужно оценить $L_{ij,k}$ – константу липшица в 2-норме градиента функции (по переменным блока $T_{(k)}$)

$$\bar f_{ij}\eta_{ij}\ln\left(\frac{1}{S+1}\left[\sum_{s\in O}\exp\left(\frac{T_{sj}-T_{si}-\bar t_{ij}}{\eta_{ij}}\right)+1\right]\right)=$$

$$=\bar f_{ij}\max_{\substack{u_0+\sum_{s\in O}u_s=1\\u_0,u_s\ge0,\ s\in O}}\left\{\sum_{s\in O}\left(T_{sj}-T_{si}-\bar t_{ij}\right)u_s-\eta_{ij}u_0\ln u_0-\eta_{ij}\sum_{s\in O}u_s\ln u_s\right\}.$$

Из выписанного представления и теоремы 1 работы [40] имеем, что

$$L_{ij,k}=\bar f_{ij}/\eta_{ij},\ k=i,j;\ L_{ij,k}=0,\ k\ne i,j.$$

Введя

$$R_T^2 = \frac{1}{2}\left\| T^0 - T_* \right\|_2^2,$$

$C_{ij} \le 10$ – число “соседей” в транспортном графе у вершин $i$ и $j$,

получим, что алгоритм 2 APPROX из работы [43] с $n$ блоками (размер каждого блока $S$), с евклидовой прокс-структурой в композитном варианте, с композитом

$$-\sum_{\substack{s\in O,\, k\in D\\ (s,k)\in W}} d_{sk}\cdot\left(T_{sk} - T_{ss}\right),$$

имеет, согласно теореме 3 [43], следующую оценку (в среднем) общей сложности решения задачи (7) с точностью по функции $\varepsilon/4$ ($\bar{L}, \bar{C}$ – специальным образом “взвешенные средние” констант $\bar{f}_{ij}/\eta_{ij}$, $C_{ij}$, причем $\bar{L} \le \max_{(i,j)\in E} \bar{f}_{ij}/\eta_{ij}$, $\bar{C} \le \max_{(i,j)\in E} C_{ij}$)

$$\underbrace{\mathrm{O}\left(n\sqrt{\frac{\bar{C}\bar{L}R_T^2}{\varepsilon}}\right)}_{\text{число итераций}}\cdot\underbrace{\mathrm{O}\left(\bar{C}S\right)}_{\substack{\text{стоимость}\\ \text{итерации}}} = \mathrm{O}\left(\bar{C}Sn\sqrt{\frac{\bar{C}\bar{L}R_T^2}{\varepsilon}}\right).$$

Если бы мы использовали не APPROX, а быстрый градиентный метод в композитном варианте [27, 40], то оценка была бы хуже

$$\mathrm{O}\left(Sn\sqrt{\frac{n\bar{L}R_T^2}{\varepsilon}}\right).$$

Заметим, что алгоритм можно распараллелить на $\mathrm{O}\left(n/\bar{C}\right)$ процессорах, каждый процессор при этом должен сделать

$$\mathrm{O}\left(\bar{C}^2 S\sqrt{\frac{\bar{C}\bar{L}R_T^2}{\varepsilon}}\right)$$

арифметических операций. При $S \ll n$ получается довольно оптимистичная оценка.

## 5. Заключение

В статье были описаны различные способы поиска равновесного распределения потоков в моделях Бэкмана (раздел 2), в модели стабильной динамики Нестерова–де Пальма (раздел 4) и промежуточных моделях (раздел 3). Промежуточные модели получаются из модели Бэкмана при предельном переходе $\mu \to 0+$ по части ребер (см. раздел 2). Модель стабильной динамики получается, когда предельный переход осуществлен на всех ребрах (см. раздел 3). Подход раздела 2 не применим к моделям разделов 3, 4, поскольку (см. раздел 2)

$$L_2 = \max_{e\in E} \tau_e'^{\mu}\left(\hat{f}_e\right) \to \infty \text{ при } \mu \to 0+.$$

Подход раздела 4 не применим к моделям Бэкмана, поскольку, по сути, ограничивается только решением специальной задачи линейного программирования. Подход раздела 3 применим ко всем рассмотренным в статье моделям.

Введя относительную точность по функции $\tilde{\varepsilon}$ (см. раздел 2), с которой мы хотим искать равновесие (в реальных приложениях относительной точности $\tilde{\varepsilon} \sim 0.01-0.05$ оказывается более чем достаточно), запишем оценки общего времени работы методы из разделов 2 – 4:

| метод / модель | **Бэкман** | **Нестеров–де Пальма** | **Промежуточная** |
|---|---|---|---|
| **Франк–Вульф** | $K_2(S,n)Sn/\tilde{\varepsilon}$ | $\varnothing$ | $\varnothing$ |
| **Зеркальный спуск** | $K_3(S,n)n/\tilde{\varepsilon}^2$ | $K_3(S,n)n/\tilde{\varepsilon}^2$ | $K_3(S,n)n/\tilde{\varepsilon}^2$ |
| **Покомпонентный спуск** | $\varnothing$ | $K_4(S,n)Sn/\tilde{\varepsilon}$ | $\varnothing$ |

**Таблица 1**

К сожалению, константы $K$ этих методов могут существенно отличаться и зависеть, в частности, от $n$. В худших случаях может быть $K = \mathrm{O}(n)$.

Заметим, что ожидаемую (среднюю) сложность $\tilde{\mathrm{O}}\left(\max\{n,|W|\}^3\right)$ поиска равновесия в модели стабильной динамики дает симплекс метод [45 – 47] (для прямой задачи – задачи линейного программирования). Из таблицы 1 можно заключить, что симплекс метод будет конкурентоспособным лишь при небольшом числе корреспонденций $|W| \le n$. Отметим при этом, что общую сложность $\tilde{\mathrm{O}}\left(n^3\right)$ имеет заметно более простая транспортная задача линейного программирования с $n$ пунктами производства и потребления, причем эта оценка является не улучшаемой [48, 49].

Интересно сравнить приведенные в таблице 1 оценки с оценками сложности решения задачи поиска стохастического равновесия в модели Бэкмана из работы [50].

В заключение сделаем существенную для практики оговорку. На полученные в данной статье оценки следует смотреть исключительно с точки зрения качественного понимания сложности того или иного метода, но ни коим образом не с точки зрения отбора лучшего метода. В данной работе отсутствуют сравнительный анализ констант методов из таблицы 1. Такой анализ требует проработки некоторых технических деталей, связанных с дополнительным погружением в специфику постановки. При выбранном в данной статье уровне грубости получения оценок можно считать, что все три метода (из разделов 2 – 4) конкурентоспособны. Помочь отобрать наилучший метод здесь могут численные экспе-рименты. Этому планируется посвятить отдельную работу.

Авторы выражают благодарность А.С. Аникину, А.С. Немировскому, Е.А. Нурмин-скому, Н.Б. Шамрай, А.А. Шананину, С.В. Шпирко за ценные замечания, а также студен-там ФУПМ МФТИ Александру Бондаренко, Марине Казначеевой, Даниилу Меркулову, Ильнуре Усмановой и студенту факультета компьютерных наук ВШЭ Антону Родоманову за помощь в работе.

Особую благодарность авторам хотелось бы выразить проф. Юрию Евгеньевичу Не-стерову, за постоянное внимание к работе и ценнейшую обратную связь.

Исследование в разделе 2 выполнено в ИППИ РАН за счет гранта Российского науч-ного фонда (проект № 14-50-00150), исследование в разделе 3 выполнено при поддержке гранта 15-31-20571-мол_а_вед, исследование А.В. Гасникова в разделе 4 выполнено при поддержке гранта РФФИ 15-31-70001 мол_а_мос.

## Приложение

Приведем с доказательство использованного нами в п. 4 результата (6) (Дорна–Нестерова, 2014).

**Теорема 5.** *Задача (1) эквивалентна задаче*

$$\min_T \left\{ -\sum_{\substack{s\in O,\, k\in D \\ (s,k)\in W}} d_{sk}\cdot\left(T_{sk}-T_{ss}\right)+\sum_{(i,j)\in E}\bar{f}_{ij}\max_{s\in O}\left(T_{sj}-T_{si}-\bar{t}_{ij},0\right)\right\}. \tag{8}$$

*Решение задачи (8) связано с решением задачи (1) следующим образом*

$$t_{ij}=\max\left\{\max_{s\in O}\left(T_{sj}-T_{si}\right),\bar{t}_{ij}\right\},\ (i,j)\in E. \tag{9}$$

**Доказательство.** Двойственная задача для задачи (1) имеет вид

$$\min_{t\geq\bar{t}}\left\{-\sum_{(s,k)\in W} d_{sk}\cdot T_{sk}\left(t\right)+\sum_{(i,j)\in E}\bar{f}_{ij}\cdot\left(t_{ij}-\bar{t}_{ij}\right)\right\}. \tag{10}$$

Введем набор переменных $T_{sk}$ – время проезда по кратчайшему пути из вершины $s$ в вершину $k$. Тогда для любых трех вершин $s$, $i$, $j$, таких, что существует ребро $(i,j)\in E$ выполнено соотношение

$$T_{sj} \le T_{si} + t_{ij} .$$

Следовательно, задачу (10) можно переписать в виде

$$\min_{t,T} \left\{ -\sum_{(s,k)\in W} d_{sk} \cdot T_{sk} + \sum_{(i,j)\in E} \bar{f}_{ij} \cdot \left( t_{ij} - \bar{t}_{ij} \right) \middle| t \ge \bar{t};\ T_{sj} \le T_{si} + t_{ij}, s \in V, (i,j) \in E \right\} . \qquad (11)$$

Это задача ЛП. В её целевом функционале все компоненты вектора $t$ неотрицательны, при этом ограничение на каждую компоненту имеет вид

$$t_{ij} \ge \max\left\{ \max_{s} \left[ T_{sj} - T_{si} \right], \bar{t}_{ij} \right\}. \qquad (12)$$

Следовательно, задача (11) может быть явно решена относительно вектора $t$ – в (12) имеет место равенство, т.е. имеет место (9). Исключая вектор $t$, приходим к задаче (8), что завершает доказательство. □

## Литература

# Numerical Methods for the Problem of Traffic Flow Equilibrium in The Beckmann and the Stable Dynamic Models

*Gasnikov A.V.*
(PreMoLab MIPT, IITP RAS)
gasnikov.av@mipt.ru
*Dvurechensky P.E.*
(WIAS, IITP RAS)
dvurechensky@iitp.ru
*Dorn Yu.V.*
(Research Centre for Transport Policy Studies, Institute for Transport Economics and Transport Policy Studies, National Research University Higher School of Economics, PreMoLab MIPT, IITP RAS)
dorn.iuv@mipt.ru
*Maximov Yu.V.*
(PreMoLab MIPT, IITP RAS)
yury.maximov@phystech.edu

**Abstract**

In this work we propose new computational methods for transportation equilibrium problems. For Beckmann's equilibrium model we consider Frank–Wolfe algorithm in a view of modern complexity results for this method. For Stable Dynamic model we propose new methods. First approach based on mirror descent scheme with euclidean prox-structure for dual problem and randomization of a sum trick. Second approach based on Nesterov's smoothing technique of dual problem in form of Dorn–Nesterov and new implementation of randomized block-component gradient descent algorithm.